\documentclass{article}
\usepackage{amsfonts}
\usepackage{amssymb}
\usepackage{latexsym}
\font\cyr=wncyr10 
\def\MW{Mordell-\mbox{\kern-.16em}Weil}
\newcommand{\nl}{\hspace*{\fill}\\}
\newcommand{\TS}{Tate-\v Safarevi\v c}
\newcommand{\BW}{Barnes-\mbox{\kern-.16em}Wall}
\newcommand{\be}{\begin{equation}}
\newcommand{\ee}{\end{equation}}
\newcommand{\bea}{\begin{eqnarray}}
\newcommand{\eea}{\end{eqnarray}}
\newcommand{\ra}{\rightarrow}
\newcommand{\Ra}{\Rightarrow}
\newcommand{\teq}{\!=\!}
\newcommand{\leqs}{\leqslant}
\newcommand{\geqs}{\geqslant}

\newcommand{\R}{\mathbf{R}}
\newcommand{\Z}{\mathbf{Z}}

\newcommand{\Ht}{\hat{h}}
\newcommand{\Sha}{\mbox{\cyr X}} 
\newcommand{\End}{\mathop{\rm End}\nolimits}
\newcommand{\Aut}{\mathop{\rm Aut}\nolimits}
\newcommand{\Hom}{\mathop{\rm Hom}\nolimits}
\newcommand{\Jac}{\mathop{\rm Jac}\nolimits}
\newcommand{\Fbar}{{\overline{\mathbf{F}_{\!2}\!}}\,}
\newcommand{\0}{^{\phantom0}}
\begin{document}

\begin{center}
\begin{Large}
\MW\ lattices in characteristic~2 \\
III: A \MW\ lattice of rank 128

\vspace*{2ex}
\end{Large}

Noam D. Elkies\\
November, 2000
\end{center}

\vspace*{4ex}

\textbf{Introduction.}  In the first paper~\cite{E1} of this series
we constructed a family of lattices in dimensions $2^{n+1}$ for
positive integers~$n$.  From the theory of elliptic curves
over function fields we obtained upper bounds on the discriminants
and lower bounds on the minimal norms of these lattices,
showing that their associated lattice packings of spheres
equal the previous records for $n \leqs 4$ and exceed them for
$5 \leqs n \leqs 9$.  We then showed, for $n=5$
(the first case of a new record), that our lower bound
on the density of the lattice packing in $\R^{64}$ is sharp,
and reported on the computation of the kissing number of the lattice,
which was at the time the largest kissing number known in $\R^{64}$.
Both of these records have since been superseded by G.~Nebe~\cite{Nebe},
and the kissing number was pushed still higher
by a nonlattice packing \cite{EdRS}.
Thus $n=6$ is now the first case in which the construction of~\cite{E1}
yields a previously unknown lattice of record density.

In this paper we analyze this $128$-dimensional lattice $MW_{128}$.
We determine its density, again showing
that the lower bound from~\cite{E1} is sharp
by proving that the elliptic curve (see (\ref{E128}) below)
has trivial \TS\ group and nonzero rational points of height~$22$,
as small as possible by \cite[Prop.2]{E1}.
We then report on a computation that determined all the rational points
of that minimal height, thus obtaining the kissing number of $MW_{128}$.
Like the packing density, the kissing number of $MW_{128}$
is by a considerable factor the largest known kissing number
{\em of a lattice} in this dimension,
though once more \cite{EdRS} gives a much larger
nonlattice kissing configuration.

$MW_{128}$ is the \MW\ lattice of an elliptic curve~(\ref{E128})
over the rational function field $k(t)$,
where $k$ is a finite field of $2^{12}$ elements.
To compute the minimal vectors of~$MW_{128}$
we listed all solutions in~$k$\/
of a system of simultaneous nonlinear equations in several variables
(the coefficients of~$x$ as a polynomial in~$t$).
We reduced the search space by using the automorphisms of $MW_{128}$
described in~\cite{E1} and solving for some of the variables.
This left about $2\cdot 10^{12}$ possibilities on which to check
the remaining, now much more complicated, equations.
An exhaustive search over this space
would still take months on a single fast computer.
Fortunately the first of these equations can be written as a quadratic
in one of the variables; this reduced the search space
by three orders of magnitude, to the point that the computer time
was comparable to the time it took to program the search.
The relatively simple form of that equation, though welcome,
was unexpected and is still unexplained.  This is one of several open
questions we raise in the concluding section of the present paper.

\pagebreak

\textbf{Statement of results.}

To simplify notation we henceforth denote the lattice $MW_{128}$ by~$M$.

Let $k$\/ be a finite field of $2^{12}\teq4096$ elements,
let $K$\/ be the rational function field $k(t)$,
and let $E/K$\/ be the potentially constant elliptic curve
\be
y^2 + y = x^3 + t^{65} + a_6,
\label{E128}
\ee
where $a_6$ is any of the $2^{11}$ elements of~$k$\/
whose absolute trace $\sum_{m=0}^{11} a_6^{2^m}$ equals~$1$.
(As noted in \cite{E1,E2},
all these choices of $a_6$ yield isomorphic curves;
if we worked over $\overline k \teq \Fbar$ instead of~$k$,
we could drop $a_6$ entirely
and simplify the equation of~$E$\/ to $y^2 + y = x^3 + t^{65}$.)
Let $M$\/ be the \MW\ group of~$E$, consisting of $\mathbf{0}$
together with solutions $(x,y)\in K \times K$\/ of~(\ref{E128}).
The map
\be
\phi: E \ra E, \
(x,y) \mapsto (x^2, y^2 + t^{65} + a_6)
\label{phi}
\ee
is an inseparable $2$-isogeny whose square is multiplication by~$-2$.
By \cite[Thm.2]{E1}, the image of~$\phi$ on $E(K)$
is the kernel of the map
\be
\epsilon: E(K) \ra K/K^2, \qquad
\mathbf{0} \mapsto 0,\ (x,y) \mapsto x \bmod K^2
\label{epsilon}
\ee
in which $K,K^2$ have the structure of additive groups.

\textbf{Theorem.}
{\em
i) $M$\/ has rank $128$ and trivial torsion.
\nl
ii)
The canonical height $\Ht$ on~$M$\/ is given by
$\Ht(\mathbf{0})=0$ and $\Ht(x,y)=\deg x$, the degree of~$x$\/
as a rational function of~$t$.  This height gives $M$\/ the structure
of an even integral lattice in Euclidean space of dimension~$128$.
\nl
iii) The Selmer group for $\phi$ is the subspace~$S$\/ of $K/K^2$
represented by the polynomials of the form
\be
(x_{21}\0 t^{21} + x_{21}^4 t^{19} + x_{21}^{16} t^{11})
+ (x_{17}\0 t^{17} + x_{17}^4 t^3)
+ x_{13} t^{13} + x_9 t^9 + x_5 t^5 + x_1 t
\label{Sel_pol}
\ee
with
\be
x_j\in k, \quad x_{13}^{16}=x_{13}\0;
\label{Sel_cond}
\ee
$S$ is an elementary abelian $2$-group of rank
\be
12+12+4+12+12+12 = 64.
\label{Sel_phi_rank}
\ee
iv) The \TS\ group of $E/K$\/ is trivial.
The discriminant of~$M$\/ is $2^{120}$.
v) The minimal norm of~$M$\/ is $22$, attained by $(x,y)$
if and only if $x,y$ are polynomials in~$t$
of degrees $22$ and~$33$ respectively.  There are
\be
218044170240 = 2^{17} \; 3 \; 5 \; 13 \; 19 \; 449
\label{kiss128}
\ee
vectors of this minimal norm in~$M$.
\nl
vi) The normalized center density of~$M$\/ is
\be
11^{64} / 2^{124} = 2^{97.4036+}.
\label{delta128}
\ee
}

\pagebreak

\textbf{Proof of parts (i) through (iv): rank, discriminant, and $\Sha$}

Each of (i), (ii), and the implication (iii)$\Ra$(iv)
is contained in the special case $(n,q)=(6,64)$
of our results in~\cite{E1}.
We briefly go over these in the next two paragraphs.

Define curves $C,E_0$ over~$k$ by
\be
C: u^2 + u = t^{65}, \qquad
E_0: Y^2 + Y = X^3 + a_6.
\label{C,E0}
\ee
Then $E_0$ is a supersingular elliptic curve, and
$C$\/ is a hyperelliptic curve of genus~$32$ whose Jacobian $\Jac(C)$
is isogenous with $E_0^{32}$~\cite[Prop.1]{E1}.
The $K$-rational points of~$E$\/ correspond bijectively
with maps from~$C$\/ to~$E_0$
that take the point at infinity of~$C$\/ to the origin of~$E_0$:
such a map is either constant or of the form
$(t,u)\mapsto(x(t),y(t)+u)$ with $(x,y)$ a nonzero point of~$E$.
This correspondence respects the group laws on~$E$\/ and~$E_0$,
and yields an identification of~$M$\/ with $\Hom(\Jac(C),E_0)$,
a group of the same rank as $\Hom(E_0^{32},E_0)=\End(E_0)^{32}$.
Thus $M$\/ has rank $4\cdot 32 = 128$, as claimed in~(i).
The formula for $\Ht$
and the fact that $\Ht(P)\in2\Z$ for all $P\in M$\/
are the case $(n,q)=(6,64)$ of \hbox{\cite[Prop.2]{E1}}.
 
The discriminant of the \MW\ lattice of an elliptic curve
over a global field is related to the order of the curve's
\TS\ group by the conjecture of Birch and Swinnerton-Dyer.
In our case of a curve over a function field, this conjecture
was formulated by Artin and Tate~\cite{Tate}
and proved under certain hypotheses by Milne~\cite{Milne}.\footnote{
  As noted in~\cite{E2}, Milne had to also assume odd characteristic,
  but this assumption was later eliminated
  by work of Illusie~\cite{Illusie},
  so we may use Milne's results also in our characteristic-$2$ setting.
  }
In~\cite[Thm.1]{E1} we observed that these hypotheses were satisfied
by each of our curves~$E$, and computed the resulting relationship
between the order of the \TS\ group $\Sha(E)$ and the discriminant 
$\Delta$ of its \MW\ lattice.  Their product is always a power of~$2$,
so $\Sha(E)$ is a $2$-group,
and is trivial if and only if $\Sha_\phi$ is trivial;
that is, if and only if $E(K)/\phi(E(K))$
is all of the Selmer group for~$\phi$.
Now $E(K)/\phi(E(K))$ is an elementary abelian \hbox{$2$-group}
whose rank is half the rank of $E(K)$, because $\phi^2=-2$.
In our present case of $q=64$, we already know that the half-rank
is $128/2=64$, so once we prove (iii) the triviality of~$\Sha(E)$
will follow.  The formula of~\cite[Thm.1]{E1} gives
\be
|\Sha(E)| \, \Delta = 2^{120},
\label{120}
\ee
for $q=64$, so the discriminant claim of~(iv) will follow as well.

It remains to verify that the $\phi$-Selmer group
is given by~(\ref{Sel_pol},\ref{Sel_cond}).
The analysis proceeds as in \cite{E1} (for $y^2+y=x^3+t^{33}$)
and \cite{E2} (for $y^2+y=x^3+t^{13}+a_6$),
but takes more steps to complete.  As happened there, it is enough
to show that $S$\/ is contained in the Selmer group,
because it has the correct size $2^{128/2}$.  The Selmer group
consists of $0$ together with all elements of $K/K^2$ that represent
the \hbox{$x$-coordinate} of a solution of (\ref{E128}) with
$x,y$ in $k((t^{-1}))$, the completion of~$K$\/ at the place $t=\infty$.

By \cite[Thm.2]{E1}, every element of the $\phi$-Selmer group
has a unique representative $\xi$ that is an odd polynomial in~$t$
(that is, a \hbox{$k$-linear} combination of $t^j$
for odd positive integers~$j$) whose degree~$d$\/ satisfies
$3d < 65$ and $d \equiv 65 \bmod 4$.  Thus $d$\/ is one of
$1,5,9,13,17,21$.  We give the proof in the case $d=21$,
which is also relevant to our computation of the minimal vectors.
The other cases are similar but easier.
Alternatively, once we find a single $P_0\in M$\/ with $\Ht(P)=22$,
and thus with $\epsilon(P_0)$ represented by a polynomial
of degree~$21$, we know that for any other point~$P$\/
at least one of $\epsilon(P)$ and $\epsilon(P+P_0)$ has $d=21$;
so once we have done $d=21$ the other cases will follow.

Suppose $x = \sum_{j=-\infty}^{d'} x_j t^j$ is the $x$-coordinate
of a point of~$E$\/ over $k((t^{-1}))$, with $x_{21}\neq 0$.
Necessarily $d'=(65-d)/2=22$.
Let
\be
\eta := x^3 + t^{65} + a_6 = \sum_{j=-\infty}^{66} \eta_j t^j.
\label{eta}
\ee
since
\be
x^3 = x \cdot x^2 =
\left( \sum_{j=-\infty}^{22} x_j t^j \right)
\left( \sum_{j=-\infty}^{22} x_j^2 t^{2j} \right),
\label{xxx}
\ee
we have $\eta_{65} = x_{21}\0 x_{22}^2 + 1$ and
\be
\eta_j = \sum_{j_1+2j_2=j} x_{j_1}\0 x_{j_2}^2
\label{etaj}
\ee
for all $j\neq0,65$.
By the Lemma in~\cite{E1}, if $\eta$ is of the form $y^2+y$
for some $y\in K((t^{-1}))$ then
\be
\sum_{m=0}^\infty (\eta\0_{2^m j_0})^{2^{-m}} = 0
\label{eta_cond}
\ee
for each odd $j_0>0$ (note that this sum is finite,
and the fractional exponents make sense in~$k$).
Taking $j_0=65,63,61,\ldots$ \ in~(\ref{eta_cond}),
and using \textsc{macsyma} to simplify the resulting equations,
we find:
\bea
\kern-2ex
j_0 = 65: && x_{22}\0 = x_{21}^{-1/2}     \label{x22} \\
\kern-2ex
j_0 = 63: && x_{19}\0 x_{22}^2 + x_{21}^3 = 0
  \ \Ra\  x_{19} = x_{21}^4     \label{x19} \\
\kern-2ex
j_0 = 61: && x_{17}\0 x_{22}^2 + x_{19}\0 x_{21}^2 + x_{21}\0 x_{20}^2 = 0
  \nonumber \\
&& \Ra x_{20}\0
     = \left( x_{21}^{-1} (x_{21}^{-1} x_{17}\0 + x_{21}^7) \right)^{1/2}
     = x_{21}^{5/2} + x_{21}^{-1} x_{17}^{1/2}
  \label{x20} \\
\kern-2ex
j_0 = 59: && x_{15}\0 = x_{21}\0 \left(
     x_{17}\0 x_{21}^2 + (x_{21}^9 + x_{21}^2 x_{17}\0) + x_{21}^9
     \right) = 0
  \label{x15} \\
\kern-2ex
j_0 = 57: && x_{18}\0 = a_{21}^{11/2} + a_{21}^2 a_{17}^{1/2}
     + a_{21}^{-3/2} a_{17}\0 + a_{21}^{-1} a_{13}^{1/2}
  \label{x18} \\
\kern-2ex
j_0 = 55: && x_{11}\0 = x_{21}\0 \left(
     x_{13}\0 x_{21}^2 + x_{17}\0 x_{21}^8 + x_{21}^4 x_{18}^2
     + x_{21}\0 x_{17}^2 \right)
     = x_{21}^{16}
  \label{x11} \\
\kern-2ex
j_0 = 53: && x_{16}\0 = a_{21}^{17/2} + a_{21}^5 a_{17}^{1/2}
     + a_{21}^{-2} a_{17}^{3/2} + a_{21}^2 a_{13}^{1/2}
     + a_{21}^{-1} a_9^{1/2}
  \label{x16} \\
\kern-2ex
j_0 = 51: && x_7\0 = x_{21}\0 \sum_{j=15}^{21} x_j^2 x_{51-2j}\0
     = \cdots = 0
  \label{x7} \\
\kern-2ex
j_0 = 49: && x_{14}\0
     = x_{21}^{-1/2} \sum_{j=15}^{22} x_j\0 x_{49-2j}^{1/2}
     = x_{21}^{23/2} + x_{21}^8 x_{17}^{1/2} + x_{21}^{9/2} x_{17}\0
     \nonumber \\
&& \kern5ex {} + x_{21}^{-5/2} x_{17}^2
     + x_{21}^{-2} x_{13}^{1/2} x_{17}\0
     + x_{21}^5 x_{13}^{1/2} + x_{21}^{-3/2} x_{13}\0
  \label{x14} \\
&& \kern5ex {} + x_{21}^2 x_9^{1/2} + x_{21}^{-1}  x_5^{1/2}
  \nonumber \\
\kern-2ex
j_0 = 47: && x_3\0 = x_{21}\0 \sum_{j=13}^{21} x_j^2 x_{47-2j}\0
     = \cdots = x_{17}^4
  \label{x3}
\eea
In particular, for each odd $j$ the coefficient $x_j$ depends on
the $\xi$~coefficients $x_1,x_5,x_9,x_{13},x_{17},x_{21}$
according to (\ref{Sel_pol}).
To finish the proof of~(iii) we must also show that $x_{13}$ is
in the \hbox{$16$-element} subfield of~$k$.  Continuing our computation
we find that the conditions (\ref{eta_cond}) with $j_0=45,41,37$
yield (increasingly complicated) formulas for $x_{12},x_{10},x_8$
in terms of the $\xi$~coefficients, while the conditions
with $j_0=43,39,35$ are satisfied automatically.  When $j_0=33$,
the first case in which the sum in~(\ref{eta_cond}) has more than
one term, we obtain a much longer expression for $x_6$.  When this
and our previous formulas are substituted into the $j_0\teq31$ condition
$\eta_{31}^2=\eta_{62}\0$, all but two of the terms cancel, leaving only
$x_{13}^{16}+x_{13}\0=0$, and we are done.

This massive cancellation and the simple equations for
$x_3,x_7,x_{11},x_{15},x_{19}$ are in striking contrast
to the increasingly complicated formulas for $x_j$ with $j$ even.
But the odd-order coefficients are constrained by the requirement
that they constitute a group.  Moreover, the Selmer group inherits
the symmetries of~$M$\/ coming from the automorphisms of~$C$\/
noted in \cite[Eqn.10]{E1}.  Namely, if $x(t)\bmod K^2$
is in the Selmer group, then so is $x(at+b) \bmod K^2$
for all $a,b\in k$\/ such that $a^{65}=1$.  This severely
constrains the possibilities for the Selmer group; for instance,
$x_{15}$ must vanish, and $x_{19}$ must be proportional to $x_{21}^4$,
else we could use linear combinations of $x(at)$ with $a^{65}=1$
to obtain a solution in $k((t^{-1}))$ of~(\ref{E128}) with $x$
a square plus a polynomial of degree $19$ or~$15$, contradicting
\hbox{\cite[Thm.2]{E1}.}  (Similar arguments arise in Dummigan's
investigation \cite{Dummigan} of the \TS\ groups of certain
constant elliptic curves related to those of~\cite{E1}.)
The even-order coefficients need not constitute a group,
but are still constrained by the invariance
under the $65\cdot 2^{12}$ transformations $t\mapsto at+b$\/;
this provides a sanity check on our formulas for those coefficients.

\textbf{Parts (v) and (vi): minimal norm, density, and kissing number}

By \cite[Thm.1]{E1}, any nonzero $(x,y)\in M$\/ has height
at least $22$, with equality if and only if $x,y$ are polynomials
in~$t$\/ of degrees~$22$ and~$33$.  Thus we can prove that
$M$\/ has minimal norm~$22$, and normalized norm density
given by~(\ref{delta128}), by finding a single such pair $(x,y)$.
To verify that the kissing number is given by~(\ref{kiss128})
we must enumerate all $(x,y)$ of that form.

To find these $(x,y)$, set $x_j=0$ for all $j<0$, write
$\eta_j$ ($0 \leqs j \leqs 66$) as polynomials in $x_0,\ldots,x_{22}$,
and solve (\ref{eta_cond}) for each odd $j_0 \leqs 65$
together with the equation
\be
\eta_0\0 = x_0^3 = y(0)^2 + y(0) + a_6
\label{eta_cond0}
\ee
for the constant coefficient~$y(0)$ of~$y$.  We have already used
the equations for $j_0\geqs 35$ to solve for all $x_j$ except
$x_0,x_2,x_4$ in terms of the $\xi$~coefficients;
the $j_0 \teq 29$ and $j_0 \teq 25$ equations give us $x_4$ and $x_2$
as well, and $j_0=21$ determines $x_0^4+x_0\0$.
(We already saw in \cite[Eqn.21]{E1} that,
due to automorphisms of~$E_0$ of order~$4$,
if $(x,y)\in M$\/ and $c^4=c$ then $x+c$ is also
the \hbox{$x$-coordinate} of a rational point of~$E$\/;
thus even once we know all $x_j$ for $j>0$
we can at most determine $x_0^4+x_0\0$, not $x_0\0$.)
Meanwhile, the $j_0 \teq 23$ equation simplifies to $x_{21}^{4095}=1$,
which only confirms that $x_{21}\in k^*$.  The remaining ten equations,
for $j_0=19$ through $j_0=1$, yield complicated polynomial equations
in the six $\xi$~coefficients.

Finding a single solution turns out to be easy: set $x_{21}=1$ and
$x_{17}=x_{13}=x_9=x_5=0$, when the equations for $j_0 \geqs 25$ give
\bea
x &=& t^{22} + t^{21} + t^{20} + t^{19} + t^{18} + t^{16} + t^{14}
 + x_1 t^{12} + t^{11}
 \nonumber \\
&& {} + (x_1 + 1) t^{10} + x_1 t^8 + x_1 t^6 + (x_1 + 1) t^4 + x_1^2 t^2
 + x_1^2 t + x_0,
\label{simple_x}
\eea
and the $j_0=19,21$ equations yield the conditions
\be
x_1^2 + x_1\0 + 1 = 0, \quad x_0^4 + x_0\0 = x_1\0
\label{simple_x_cond}
\ee
on the unknown coefficients $x_0,x_1$.  We calculate that
the equations (\ref{eta_cond}) for the remaining~$j_0$
are then satisfied automatically, leaving only (\ref{eta_cond0}),
which as expected has solutions $y(0)\in k$.
(In all $8$ solutions of~(\ref{simple_x_cond}),
$x_1$ and $x_0$ are of degrees $2,4$ respectively over $\mathbf{F}_2$;
thus $x_0^3$ is a fifth root of unity,
so its trace as an element of~$k$\/ equals~$1$,
whence $x_0^3+a_6\0$ is of the form $y(0)^2+y(0)$ with $y(0)\in k$.)
This proves that $M$\/ has minimal norm~$22$
and normalized center density $11^{64}/2^{124}$.

Enumerating all the minimal vectors is a more demanding task.
There are $2^{64}-2^{52}$ possibilities for the $\xi$~coefficients,
far too many for an exhaustive search.  But our equations have
many automorphisms, coming from the symmetries of~$M$\/
described in~\cite{E1} (after the proof of Prop.2).
We saw already the $65\cdot 2^{12}$ maps
\be
x(t) \mapsto x(at+b)\quad (a,b\in k,\ a^{65}=1).
\label{at+b}
\ee
We may augment these by $x\mapsto\alpha x$ for $\alpha^3$
(again inherited from $\Aut(E_0)$)
and by the twelve field automorphisms of~$k$ \cite[Eqn.22]{E1}.
This lowers the total to
$(63\cdot 2^{40})/(3\cdot 12) \approx 2 \cdot 10^{12}$.
The \hbox{$8$-element} quaternion group acting on~$E_0$
also acts on the minimal vectors,
but only by permuting the $8$ choices of $x_0$ and $y(0)$
associated to each valid sextuple of $\xi$~cofficients.
Thus these automorphisms do not further cut down our search space,
though the condition that $x_0$ and $y(0)$ must be in~$k$\/
will somewhat reduce the average work per candidate.

We might now organize the search as follows.
Initialize various tables for arithmetic in~$k$, such as
multiplication and multiplicative inverse, squaring, and exponentiation.
(For addition we use the bitwise \hbox{exclusive-or} operator
that is already built into the programming language~C.)
Using the automorphisms $t\mapsto at$, we may assume that $x_{21}$
is in the quadratic subfield $\mathbf{F}_{64}$ of~$k$.
Thanks to $x\mapsto \alpha x$, we may further limit attention
to one representative of each of the $21$ cosets
of the cube roots of unity in $\mathbf{F}_{64}^*$.
Each of these $21$ choices of $x_{21}$ then represents
$3(2^6+1)=195$ choices of~$x_{21}$ in~$k^*$.
Then, since $t\mapsto t+b$\/ fixes $x_{21}$ and translates $x_{17}$ by
\be
b^4 x_{21} + b^2 x_{19} = b^4 x_{21}\0 + b^2 x_{21}^4
= x_{21}^7 \left( (b^2 x_{21}^{-3})^2 + b^2 x_{21}^{-3} \right),
\label{tr17}
\ee
we may assume that $x_{17}$ is either~$0$ or $x_{21}^7 a_6\0$,
each possibility representing $2^{11}$ choices of $x_{17}\in k$\/
and still invariant under one translation
$t \leftrightarrow t+x_{21}^{3/2}$.
This translation does not affect $x_{13}$,
but (for most choices of $x_{21}$ and $x_{13}$) does move $x_9$.
Using this translation as well as the $12$ field automorphisms
reduces the $21 \cdot 2 \cdot 2^4 \cdot 2^{12} = 21 \cdot 2^{17}$
possibilities for $(x_{12},x_{17},x_{13},x_9)$ to slightly over
$21 \cdot 2^{17} / 24$ or about $10^5$ choices.  For each of these,
we use (\ref{x22}ff.), or the recursion (\ref{eta_cond}), to compute
$x_{22},x_{20},x_{19},x_{18},x_{16},x_{11},x_3$.
Then loop over $2^{24}$ choices of $x_5,x_1\in k$.  For each one,
solve the equations (\ref{eta_cond}) with $j_0=45,41,37,33,29,25,21$
to obtain $x_{12},x_{10},x_8,x_6,x_4,x_2$, and $x_0^4+x_0\0$.
Look up a precomputed table to choose $x_0$
and solve (\ref{eta_cond0}) for $y(0)$, if solutions exist in~$k$\/
(if not, proceed to the next $(x_5,x_1)$ pair).
If one of the eight possible $(x_0,y(0))$ is defined over~$k$,
then all are, but we need only try one because they constitute an orbit
under the \hbox{$8$-element} quaternion group.  Using the chosen $x_0$,
check whether the conditions (\ref{eta_cond}) for odd $j_0\leqs 19$
are all satisfied.  The kissing number is the sum of the orbit sizes
of the resulting minimal points of~$E$.

Unfortunately the size $2\cdot 10^{12}$ of the search space 
is too large for this computation to conclude in reasonable time.
Fortunately the $j_0=19$ equation, expanded as a polynomial
in the $\xi$~coefficients, is a quadratic equation in~$x_1^{16}$, namely
\be
A x_1^{32} + A^2 x_1^{16} = B,
\quad{\rm where}\quad
A = x_{21}^{288} x_{13}\0 + x_{21}^{192} x_9^{16}
   + x_{21}^{96} x_5^{16} + x_{21}^{496}
\label{quadratic}
\ee
and $B$\/ is
\bea
&
\kern-10ex
(x_{21}^{272} + x_{13}\0 x_{21}^{64}) x_9^{64} + x_{21}^{576} x_9^{48}
+ (x_{21}^{480} x_5^{16} + x_{21}^{880} + x_{13}\0 x_{21}^{672}
+ x_{21}^{100} + x_{13}^4 x_{21}^{48}) x_9^{32}
\kern-10ex
&
\nonumber
\\
&
+ (x_{21}^{384} x_5^{32} + x_{21}^{1184} + x_{13}^2 x_{21}^{768}) x_9^{16}
+ x_{21}^{288} x_5^{48} + (x_{21}^{688} + x_{13}\0 x_{21}^{480}) x_5^{32}
&
\nonumber
\\
&
+ (x_{13}^2 x_{21}^{672} + x_{13}^4 x_{21}^{256}) x_5^{16}
+ (x_{21}^{608} + x_{13}^2 x_{21}^{192}) x_5^4
&
\nonumber
\\
&
+ x_{13}^3 x_{21}^{864} + x_{21}^{708} + x_{13}^5 x_{21}^{448}
+ x_{13}^2 x_{21}^{292} + x_{13}^2 x_{21}^{97}.
&
\label{B_lecch}
\eea
Ugly as this may look, it is much better than looping over~$x_1$
--- especially since most of the computation of~$B$\/ can be done
independently of~$x_5$.  Except in the rare case that $A=B=0$,
we are thus left with either $2$ or~$0$ choices of~$x_1$ for each
$x_{21},x_{17},x_{13},x_9,x_5$.  It is then feasible to test every
possible $(x_{21},x_{17},x_{13},x_9,x_5,x_1)$
as outlined in the preceding paragraph.

We ran this computation and found a total of
$2940$ orbits of minimal points, with stabilizers of orders
distributed as follows:
\be
\begin{array}{c|cccccccc}
|{\rm Stab}| &    1 &   2 &  3 &  4 & 6 & 8 & 12 & 24 \\ \hline
\#           & 2766 & 134 & 21 & 11 & 3 & 1 &  3 &  1
\end{array}
\label{counts}
\ee
(For instance, (\ref{simple_x},\ref{simple_x_cond})
is in one of the orbits with a \hbox{$6$-element}
stabilizer, coming from Gal($k/\mathbf{F}_4$).
See \textsf{http://www.math.harvard.edu/\~{}elkies/mv128.txt}
for a full list of orbit representatives.)
Summing $1/|{\rm Stab}|$ over the orbits,
and multiplying the resulting total of $8531/3$
by the number $2^{12} \cdot 65 \cdot 24 \cdot 12$
of known automorphisms of~$M$, we obtain
the kissing number (\ref{kiss128}) of~$M$.~~$\Box$

\textbf{Remarks and questions.}
The unique orbit with a $24$-element stabilizer is represented
by a point~$P$\/ with $\xi$~coefficients
\be
(x_{21},x_{17},x_{13},x_9,x_5,x_1) =
(1,a,1,a+1,a^3+a^2+a,a^2+1),
\label{stab24}
\ee
where $a$ is a fifth root of unity.
This formula clearly shows three stabilizing automorphisms
(from Gal($k/\mathbf{F}_{16}$)), but in fact $P$\/ is stabilized
by two automorphisms for each element of Gal($k/\mathbf{F}_2$):
translating $t$ by either $a^3$ or $a^3+1$ has the same effect
on its \hbox{$x$-coordinate} as applying the Galois automorphism
$c\mapsto c^2$ of~$k$.  This yields a cyclic stabilizer of order~$24$.

Can the simple form $A x_1^{32} + A^2 x_1^{16} = B$\/
of~(\ref{quadratic}) be explained conceptually?
Can the kissing number (\ref{kiss128}), and more generally
the kissing numbers of our \MW\ lattices in dimensions $2^{n+1}$,
be obtained without such a long computation (which seems out of the
question already for the next case $n=7$)?  Must there always
be some nonzero vectors in the narrow \MW\ lattice whose norm
attains the lower bound $2\lfloor(2^n+4)/6\rfloor$ of~\cite[Thm.1]{E1}?
Finally, can it be shown that the $48n(q^3+q^2)$ known automorphisms
of~$MW_{2^{n+1}}$ constitute its full automorphism group once $n\geqs4$?

A final remark: multiplication of~$t$\/ by a fifth root of unity
generates an automorphism of~$M$\/ of order~$5$; the sublattice
fixed by this automorphism is the \MW\ lattice of
$y^2+y=x^3+t^{13}+a_6$, and is thus homothetic with
the Leech lattice by~\cite{E2}.

\textbf{Acknowledgements.}
This work was made possible in part by funding from
the National Science Foundation and the Packard Foundation.
I thank the Mathematical Sciences Research Institute for its
hospitality while I wrote this paper.

\vspace*{5ex}
\begin{small}
\noindent 
Dept.\ of Mathematics\\
Harvard University\\
Cambridge, MA 02138 USA\\ \\
\textsf{elkies@math.harvard.edu}
\end{small}

\end{document}